%% file: main.tex
\documentclass[numbered]{trbunofficial}
\usepackage{graphicx}
\usepackage[short]{optidef}
\usepackage[fleqn]{amsmath}
\usepackage{amssymb,amsthm}
\usepackage[ruled,linesnumbered,noresetcount,vlined]{algorithm2e}

\SetKwFor{Repeat}{repeat}{:}{endw}
\SetKwFor{Forall}{forall}{}{endw}

\usepackage[hidelinks]{hyperref}

\usepackage[colorinlistoftodos]{todonotes}

\AuthorHeaders{Lu, Riley, Gurumurthy, and Van Hentenryck}
\title{Revitalizing Public Transit in Low Ridership Areas: An Exploration of On-Demand Multimodal Transit Systems}

\author{%
  \textbf{Jiawei Lu}\\
  H. Milton Stewart School of Industrial and Systems Engineering\\
  Georgia Institute of Technology\\
  Email: jlu486@gatech.edu\\
  \hfill\break
  \textbf{Connor Riley}\\
  H. Milton Stewart School of Industrial and Systems Engineering\\
  Georgia Institute of Technology\\
  Email: ctriley@gatech.edu\\
  \hfill\break
  \textbf{Krishna Murthy Gurumurthy}\\
  Transportation and Power Systems Division\\
  Argonne National Laboratory\\
  Email: kgurumurthy@anl.gov\\
  \hfill\break
  \textbf{Pascal Van Hentenryck}\\
  H. Milton Stewart School of Industrial and Systems Engineering\\
  Georgia Institute of Technology\\
  Email: pvh@gatech.edu
  }

\TotalWords{4940}

\begin{document}
\maketitle

\section{Abstract}
Public transit plays an essential role in mitigating traffic congestion, reducing emissions, and enhancing travel accessibility and equity. One of the critical challenges in designing public transit systems is distributing finite service supplies temporally and spatially to accommodate time-varying and space-heterogeneous travel demands. Particularly, for regions with low or scattered ridership, there is a dilemma in designing traditional transit lines and corresponding service frequencies. Dense transit lines and high service frequency increase operation costs, while sparse transit lines and low service frequency result in poor accessibility and long passenger waiting time. In the coming era of Mobility-as-a-Service, the aforementioned challenge is expected to be addressed by on-demand services. In this study, we design an On-Demand Multimodel Transit System (ODMTS) for regions with low or scattered travel demands, in which some low-ridership bus lines are replaced with flexible on-demand ride-sharing shuttles. In the proposed ODMTS, riders within service regions can request shuttles to finish their trips or to connect to fixed-route services such as bus, metro, and light rail. Leveraging the integrated transportation system modeling platform, POLARIS, a simulation-based case study is conducted to assess the effectiveness of this system in Austin, Texas.

\hfill\break%
\noindent\textit{Keywords}: On-Demand Multimodel Transit System, Mobility-as-a-Service, Ride-Sharing, Transportation System Simulation
\newpage

\input{sections/1_introduction}
\input{sections/2_system_design}
\input{sections/3_methodology}
\input{sections/4_polaris}

\input{sections/5_experiments}
\input{sections/6_conclusions}
\input{sections/acknowledgements}
\newpage
\bibliographystyle{trb}
\bibliography{reference}
\end{document}

%% file: sections/1_introduction.tex
\section{Introduction}
\label{sect:introduction}

Public transit plays an indispensable role in alleviating traffic congestion, minimizing emissions, and augmenting travel accessibility and equity. Traditional public transit systems predominantly consist of fixed-route services, such as buses, and light and heavy rail, which operate on regular schedules. Over the years, researchers have extensively explored the design and operation of public transit systems, delving into areas like service line planning \cite{schobel2012line,lyu2019cb,borndorfer2007column,silman1974planning,steven2003optimization}, timetabling \cite{ibarra2016multiperiod,ibarra2012synchronization,shang2019bus,fouilhoux2016valid}, crew scheduling \cite{boyer2018vehicle,rodrigues2006vehicle,smith1988bus,ma2017fairness}, and demand estimation and prediction \cite{sun2021estimating,li2020graph}. One of the critical challenges in the design of public transit systems is the  distribution of finite service resources across time and space to meet the fluctuating and diverse travel demands. This challenge is particularly pronounced for fixed-route services, where routes remain largely unchanged after the initial planning stage.

Beyond offering a cost-efficient and eco-friendly solution for high-volume transit in busy areas, public transit systems also boost mobility and accessibility in remote areas with lower ridership. To strike a balance between operational costs and travel convenience, planners often design sparse transit lines with lower service frequency in these areas. This approach, however, may require riders to walk long distances to the nearest transit stop and endure long wait times due to the infrequent service. Such inconveniences often deter riders from using public transit, leading to decreased ridership and perpetuating a vicious cycle \cite{taylor2009nature}.

In recent years, the advent of the Mobility-as-a-Service (MaaS) has created increasing interest in on-Demand Multimodal Transit Systems (ODMTS) from both academic and industry circles. This is primarily due to ODMTS's potential to sculpt a more efficient, accessible, and equitable transportation system. Within the ODMTS framework, dynamic on-demand ride systems seamlessly integrate with conventional fixed-route public transit services. Here, traditional fixed-route systems cater to large corridor volumes, while on-demand ride systems bridge the first-mile and last-mile travel gap, thereby maximizing the benefits of both paradigms \cite{van2019demand}. Both academic research and real-world pilot projects have showcased the efficacy of ODMTS in enhancing travel experiences \cite{dalmeijer2020transfer,basciftci2023capturing,guan2023path,riley2020real}.

In this study, we explore the revitalization of public transit systems in regions with low ridership within the ODMTS framework, aiming to boost passenger travel experiences and overall system efficiency. Our primary strategy involves the introduction of dedicated on-demand ride services in these regions to connect passengers with cross-regional fixed-route transit systems, while simultaneously phasing out local low-ridership transit lines. We present an overview of the system framework and key methodologies and using the integrated transportation system modeling platform POLARIS \cite{auld2016polaris}, we design a simulation-based case study to examine the effectiveness of the proposed system in Austin, Texas. This research aspires to aid public transit design in the forthcoming era of ODMTS, particularly in regions with low traditional transit ridership.

%% file: sections/2_system_design.tex
\section{System Design}
\label{sect:system_design}

\subsection{On-Demand MultiModal Transit Systems}

Broadly speaking, public transportation can be categorized into two main mode types: conventional fixed-route transit (e.g. bus and light rail), and on-demand ride services. Each mode type has its own advantages and disadvantages. In general, the former is more cost-effective but with lower accessibility, while the latter offers more flexibility and convenience, though it tends to be more expensive and less environmentally efficient due to the individualized service. In practice, these two mode types are typically operated independently by distinct agencies, often lacking effective coordination between them.

\begin{figure}[!ht]
    \centering
    \includegraphics[width=\textwidth]{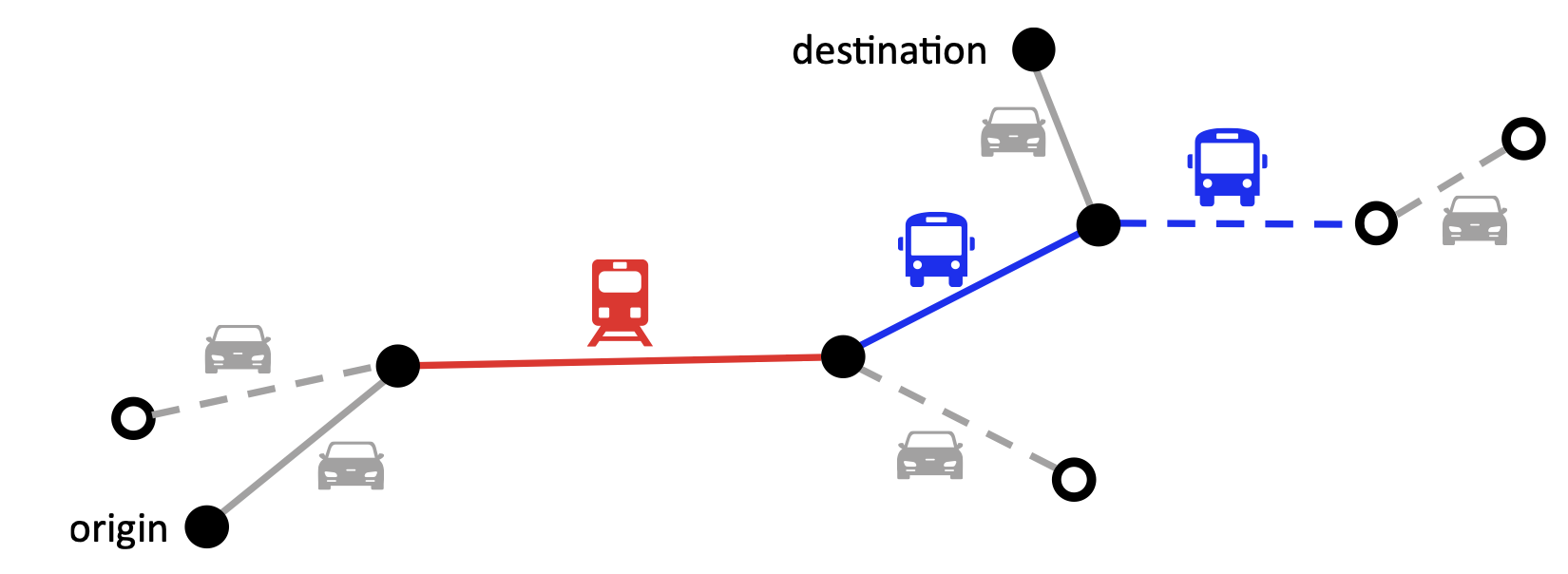}
    \caption{Example ODMTS with a passenger path (solid lines)}
    \label{fig:ex_ODMTS}
\end{figure}

On-Demand Multimodal Transit Systems (ODMTS) provide a comprehensive solution for systematically integrating existing public transportation modes, with the capability of fully taking advantage of each individual mode.
The key idea behind OMDTS is that fixed-route transit are used to serve large volumes along commute corridors while flexible on-demand ride systems are used to connect regions that are not covered by the fixed-route service to increase travel mobility and accessibility. When passengers request long-distance trips within ODMTS, routes involving potential transfers between various transit modes are automatically planned.
These plans consider factors such as cost, convenience, and the overall impact on the transportation system. For a specific passenger, on-demand shuttles assist in covering the first and last mile of the journey, thereby connecting the passenger to fixed-route transit. Once connected, passengers complete the majority of their trips in a cost-efficient manner using fixed routes. In ODMTS, as fixed routes serve aggregated large volumes, high service rates are expected to reduce passenger waiting and transferring times while ride-sharing dispatching strategies can be designed to fully utilize the capacity of shuttles and improve the efficiency of on-demand ride systems. Figure~\ref{fig:ex_ODMTS} presents an example of ODMTS with a passenger path in solid lines. 

\subsection{Public Transit System Design in Low-Ridership Regions}

In this paper, we aim to explore the potential for addressing the challenges intrinsic to traditional public transit system planning and operation in low-ridership regions, under the framework of ODMTS. These challenges originate from the trade-off between accessibility and efficiency in public transit systems. In regions with low ridership, increasing the density of traditional fixed-route transit lines and service frequency can enhance accessibility, but may result in significantly decreased cost-effectiveness and resource utilization, and vice versa. In practice, constrained by budget, there is typically a limited number of fixed-route transit lines with low service frequencies in low-ridership regions, which results in poor accessibility, long waiting times, and a low level of service.

\begin{figure}[!ht]
    \centering
    \includegraphics[width=.8\textwidth]{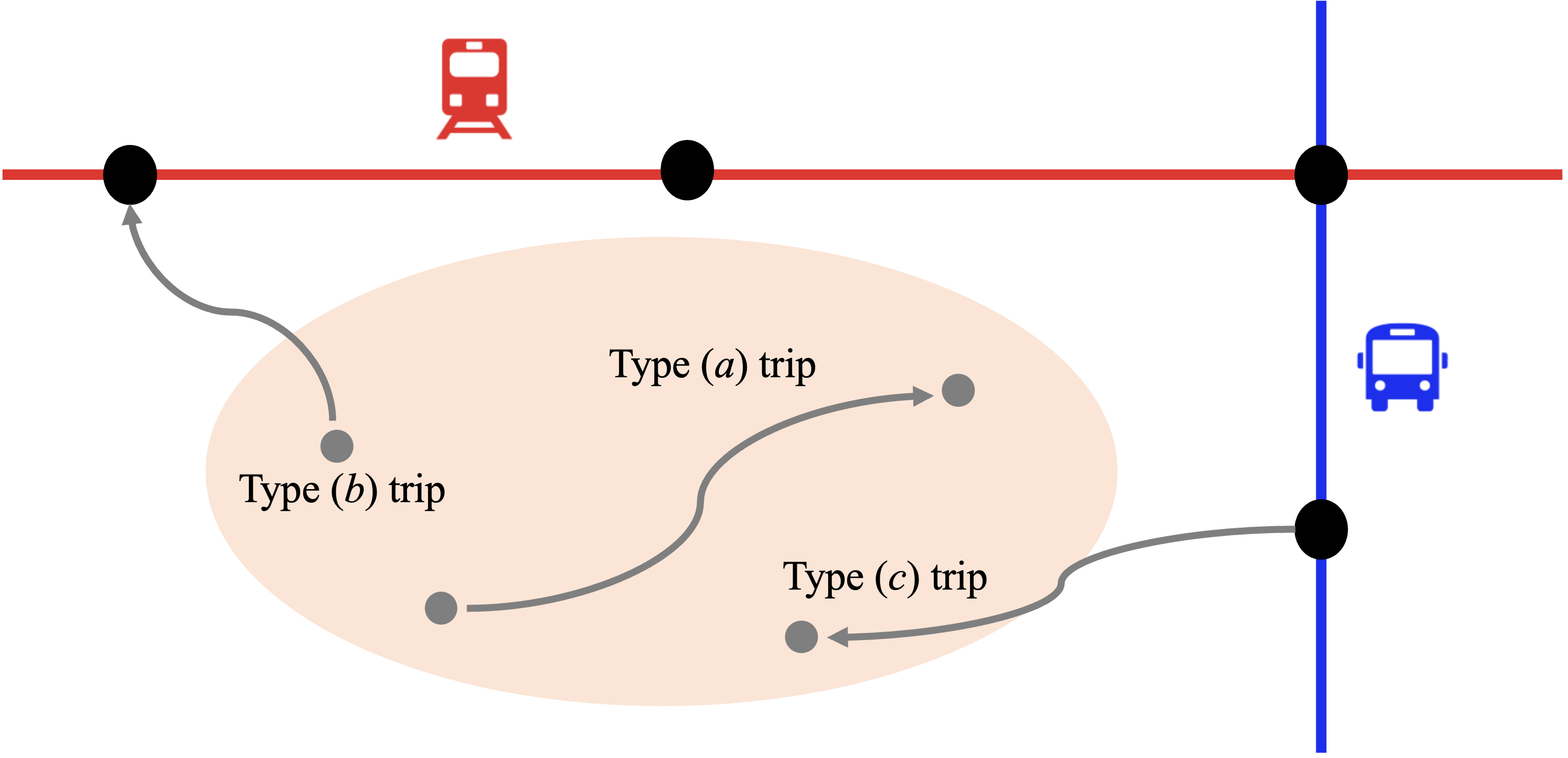}
    \caption{Three types of trips fulfilled by on-demand services in low-ridership regions}
    \label{fig:ex_geofence_trip}
\end{figure}

The key idea of our proposed design for public transit systems in low-ridership regions is substituting traditional fixed-route transit lines with dynamic, on-demand ride services. This is motivated by the observed spatial and temporal sparsity of trip requests in these regions, as well as the inherent flexibility of on-demand ride services in accommodating such requests. In the proposed system, a dedicated shuttle fleet is assigned to each low-ridership area, exclusively serving trip requests associated with that specific region. As depicted in Figure~\ref{fig:ex_geofence_trip}, within a given region, this dedicated shuttle fleet caters to three types of trips: type (\textit{a}) - trips that start and end within the region; type (\textit{b}) - trips that originate in the region and terminate at fixed-route transit stations outside the region; and type (\textit{c}) - trips that begin at fixed-route transit stations outside the region and conclude within the region. In on-demand ride services, origins and destinations of type \textit{a} trips are very straightforward, which can be the same with the pick-up and drop-off locations proposed by passengers. Types \textit{b} and \textit{c} trips are typically observed when passengers submit requests for long, cross-region journeys in ODMTS, wherein the on-demand ride service acts as a bridge connecting passengers to fixed-route transit systems. Accordingly, fixed-route transit stations outside the region are determined by the routing engine within the ODMTS. Importantly, these stations are not necessarily the closest to the original pick-up or drop-off locations specified by passengers, but rather those which facilitate the most efficient travel experience over the entirety of the trip. Given trip requests with specified pick-up and drop-off locations, the following section will illustrate how to efficiently dispatch shuttle fleets to optimize service for these requests.

%% file: sections/3_methodology.tex
\section{On-Demand Shuttle Fleet Dispatching}
\label{sect:methodology}

\newcommand{\tpmod}[1]{{\@displayfalse\pmod{#1}}}
\let\oldnl\nl
\newcommand{\nonl}{\renewcommand{\nl}{\let\nl\oldnl}}
\def\funname{\textsc{Release}}

The core of on-demand shuttle fleet dispatching lies in aptly matching trip requests with shuttle fleets and identifying optimal pick-up and drop-off sequences, with an aim to minimize overall costs. These costs may encompass shuttle fleet operational expenses, passenger travel costs, or a combination of both. In this study, our goal is to minimize passenger waiting time while maximizing the number of trips serviced.

Figure~\ref{fig:dispatching_model} presents the on-demand shuttle fleet dispatching model adopted in this study. In practical applications, this model is solved at regular intervals (e.g., every 30 seconds) to produce optimal matching of shuttles with passengers. The objective function, denoted by equation~\eqref{eq:dispatching_model_obj}, is comprised of two parts: the cost of missing a trip request and the cost associated with passenger waiting time. In the first part, the index $r \in R$ corresponds to individual trip requests, with $R$ representing the set of all trip requests within the current dispatching period. Importantly, $R$ may encompass trips from previous intervals that were not successfully matched and remain in the queue for service. The parameter $m_r$ refers to the cost incurred from missing trip request $r$, while $z_r$ is a binary variable indicating whether request $r$ was missed in the current interval. In the second part, the index $i \in \Omega$ corresponds to shuttle-trip request assignment plans, with $\Omega$ indicating the set of all such plans. Each assignment plan consists of one shuttle and its assigned trip requests. The parameter $c_i$ quantifies the passenger waiting cost for assignment plan $i$, and $y_i$ is a binary variable signifying the selection of assignment plan $i$. Constraint~\eqref{eq:request} dictates that a request $r$ can either be unfulfilled or assigned to a particular plan (i.e., a specific shuttle), where $a_r^i$ is a mapping coefficient indicating whether assignment plan $i$ services request $r$ (1) or not (0). Constraint~\eqref{eq:vehicle} ensures that a single shuttle can choose at most one assignment plan, where $v$ is the index of shuttles, $V$ is the set of shuttles, $\Omega_v$ denotes the set of assignment plans that use shuttle $v$. Finally, Constraints~\eqref{eq:z_range}~and~\eqref{eq:y_range} define the domains of the variables.

\begin{figure}[!ht]
    \begin{mini!}
	{}
        {\sum_{r \in R} m_r z_r + \sum_{i \in \Omega} c_i y_i \label{eq:dispatching_model_obj}}
	{\label{formulation:dispatching_model}}
	{}
	\addConstraint
	{\sum_{i \in \Omega} y_i a_r^i + z_r}
	{= 1 \quad \label{eq:request}}
	{\forall r \in R}
	\addConstraint
	{\sum_{i \in \Omega_v} y_i}
	{= 1 \quad \label{eq:vehicle}}
	{\forall v \in V}
	\addConstraint
	{z_r}
        {\in \{0,1\} \quad \label{eq:z_range}}
	{\forall r \in R}
	\addConstraint
	{y_i}
        {\in \{0,1\} \quad \label{eq:y_range}}
	{\forall i \in \Omega}
    \end{mini!}%
    \caption{Shuttle fleet dispatching model}
    \label{fig:dispatching_model}
\end{figure}

Model~\ref{formulation:dispatching_model} employs a column-based formulation. A key challenge in the development and resolution of Model~\ref{formulation:dispatching_model} lies in the preparation of the assignment plan set $\Omega$, given that the size of $\Omega$ escalates exponentially in correlation with the number of trip requests. Typically, scholars in the field have employed column generation as an iterative strategy for adding promising elements to set $\Omega$, avoiding the need for complete enumeration. For a more detailed exploration of column generation applied to similar problem-solving, readers are referred to \citep{riley2019column}. In this study, our focus is on shuttle fleet dispatching in independent regions with low ridership, and Model~\ref{formulation:dispatching_model} is solved at frequent, regular intervals. As a result, the number of trip requests within each dispatching interval should be constrained, and the assignment plans in set $\Omega$ can be enumerated with relative ease. Eq.\eqref{eq:assignment_set} defines set $\Omega$, where $K$ is a positive integer that controls the size of set $\Omega$ by capping the maximum number of ride-sharing trip requests in a single dispatching period. Importantly, in Eq.\eqref{eq:assignment_set}, $|G|$ is permitted to be 0, accommodating the scenario in which vehicles are not assigned new trip requests within the current dispatching period.

\begin{linenomath}
  \begin{equation}\label{eq:assignment_set}
  \Omega = \{ (v, G) \mid v \in V, G \subseteq R, 0 \leq |G| \leq K \}
  \end{equation}
\end{linenomath}

Upon establishing the assignment plan set $\Omega$, the subsequent step entails specifying the cost of each assignment plan $i \in \Omega$. In this study, we adopt passenger waiting time as the cost metric for assignment plans. Algorithm~\ref{alg:apc} elucidates the calculation of cost for a specific assignment plan within a branch-and-bound framework. 

\begin{algorithm}[!ht]
    \caption{\textsc{AssignmentPlanCost}}
    \label{alg:apc}
    \setcounter{AlgoLine}{0}

    \DontPrintSemicolon
    $cost \gets$ {\sc EvaluateAssignmentPlan}($v$, $G$) - {\sc EvaluateAssignmentPlan}($v$, $\{\}$) \;
    \SetKwFunction{FEA}{{\sc EvaluateAssignmentPlan}}
    \SetKwFunction{FRN}{{\sc CreateRootTravelSearchNode}}
    \SetKwFunction{FNS}{{\sc GetPossibleNextStops}}
    \SetKwFunction{FEN}{{\sc ExtendTravelSearchNode}}
    \SetKwProg{Fn}{Function}{:}{}
    
    \label{a1:fundecl}
    \nonl \Fn{\FEA{v, G}}
          {
            $node \gets$ {\sc CreateRootTravelSearchNode}($v$, $G$) \;
            $q \gets \{node\}$ \;
            $bestnode \gets null$, $bestcost \gets inf$ \;
            \While{$q \neq \emptyset$}
            {
                $node \gets$ select one node from $q$, remove it from $q$ \;
                $S \gets$ {\sc GetPossibleNextStops}($node$) \;
                \If{$S \neq \emptyset$}
                {
                    \Forall{$s \in S$}
                    {
                        $childnode \gets$ {\sc ExtendTravelSearchNode}($node$, $s$) \;
                        \If {$childnode.w < bestcost$}
                        {
                            $q \gets q \cup \{childnode\}$ \;
                        }
                    }
                }
                \Else
                {
                    \If{$node.w < bestcost$}
                    {
                        $bestnode \gets node$, $bestcost \gets node.w$ \;
                    }
                }
            }
            \Return $bestcost$ \;
          }
    \nonl \Fn{\FRN{v, G}}
          {
             $node \gets \{s(v), t(v), R_p(v) \cup G, R_d(v), 0\}$ \;
             \Return $node$ \;
          }
    \nonl \Fn{\FNS{node}}
          {
             $S \gets \{P(r) \mid r \in node.R_p\} \cup \{D(r) \mid r \in node.R_d\}$ \;
             \Return $S$ \;
          }
    \nonl \Fn{\FEN{node, s}}
          {
             $t \gets node.t + TravelTime(node.s,s)$ \;
             \Forall{$r \in node.R_p$}
             {
               \If{P(r) == s}
               {
                 $R_p \gets node.R_p \setminus \{r\}$,  $R_d \gets node.R_d \cup \{r\}$, $w \gets node.w + (t - RequestTime(r))$ \;
               }
             }
             \Forall{$r \in node.R_d$}
             {
               \If{D(r) == s}
               {
                 $R_p \gets node.R_p$,  $R_d \gets node.R_d \setminus \{r\}$, $w \gets node.w$ \;
               }
             }
             $childnode \gets \{s, t, R_p, R_d, w\}$ \;
             \Return $childnode$ \;
          }  
\end{algorithm}

For a given assignment plan $(v, G)$, Algorithm~\ref{alg:apc} computes the passenger waiting time that results from assigning trip requests $G$ to shuttle $v$. A travel search node in this context is defined as $node=\{s,t,R_p,R_d,w\}$, where $s$ and $t$ signify that the shuttle arrives at pick-up/drop-off stop $s$ at time $t$. Here, $R_p$ and $R_d$ correspond to trip requests awaiting pick-up and drop-off respectively, while $w$ represents the accumulated passenger waiting time. The generation of set $S$ prior to the extension of a travel search node (as conducted in line 7) is critical due to the potential for different trip requests to share the same pick-up/drop-off stops. In line 17, $s(v)$ and $t(v)$ denote the stop towards which shuttle $v$ is directed and its arrival time respectively, while $R_p(v)$ and $R_d(v)$ signify the trip requests scheduled for pick-up and drop-off by shuttle $v$ before the assignment of trip requests $G$. In essence, Algorithm~\ref{alg:apc} identifies the optimal pick-up and drop-off sequence to minimize passenger waiting time if trip requests $G$ are assigned to shuttle $v$, taking into consideration the trip requests already allocated to the shuttle. It is important to note that Algorithm~\ref{alg:apc} promotes ride-sharing to fully capitalize on the capacity of shuttles, thereby enhancing service efficiency.

%% file: sections/4_polaris.tex
\section{POLARIS: An Integrated Transportation System Modeling Platform}
\label{sect:polaris}

The on-demand dispatching algorithm was incorporated into a simulator to quantify the benefits of using such fleets in a large-scale setting. In this work, POLARIS \cite{auld2016polaris} is chosen to test the dispatching algorithm thanks to its computational speed and low runtime. POLARIS is a large-scale agent-based activity-based simulation tool that models every aspect of the transportation system as agents. This allows modelers to incorporate behavioral models for travelers that make choices when traversing the network, and add algorithms and technology to vehicles, fleet operators, and other aspects of the supply to capture its impact appropriately. 

\begin{figure}[!ht]
    \centering
    \includegraphics[width=\textwidth]{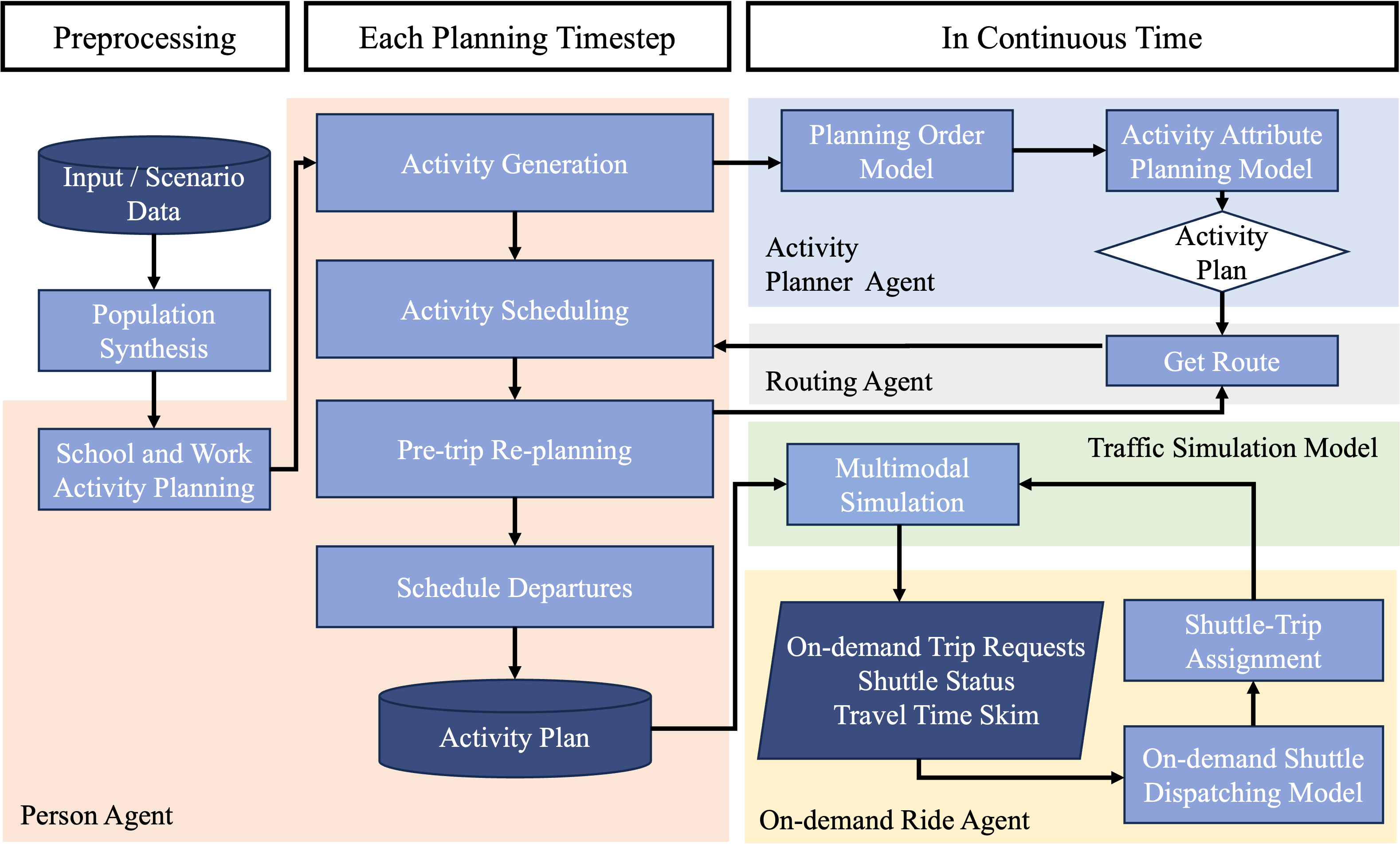}
    \caption{Embedded on-demand shuttle dispatching in POLARIS framework, adapted from \cite{auld2016polaris}}
    \label{fig:polairs}
\end{figure}

POLARIS can be generalized at a high level to comprise of the demand and supply aspects. For a given simulation run, a population is synthesized based on the underlying demographic distribution in a region (provided by Census, ACS, and PUMS), and econometric models are evaluated to generated activities, set their duration, mode, destination, and a conflict resolving algorithm is used to generate the demand for simulation \cite{auld2009framework}. A Newells traffic flow model \cite{de2019mesoscopic} and an intermodal A* router \cite{verbas2018time} provide the necessary algorithms to assign the demand on the network. Additional modules to model transit, transportation network company \cite{gurumurthy2020integrating}, and freight \cite{stinson2022introducing} are included to gather a comprehensive understanding of travel in the system. High-performance C++ code is used in POLARIS which allows for a large region with its entire population to be simulated in under 2 hr per run.  
The shared mobility module in POLARIS \cite{gurumurthy2020integrating} was updated to include the dispatching algorithm outlined in the previous section as a new strategy. This in turn provides the ability to benchmark and compare against existing dispatching algorithms used in POLARIS. In order to isolate the randomness that arises from synthesizing a population and generating its demand, all demand can be fixed after an initial POLARIS run for a region. Fixed demand generates the necessary background congestion to identify the change in impact from solely the dispatching algorithm that is used here. Figure~\ref{fig:polairs} presents the interactions between the proposed on-demand shuttle dispatching model with other modules in POLARIS.

%% file: sections/5_experiments.tex
\section{Case Study}
\label{sect:experiment}

\begin{figure}[!ht]
\centering
\includegraphics[width=\textwidth]{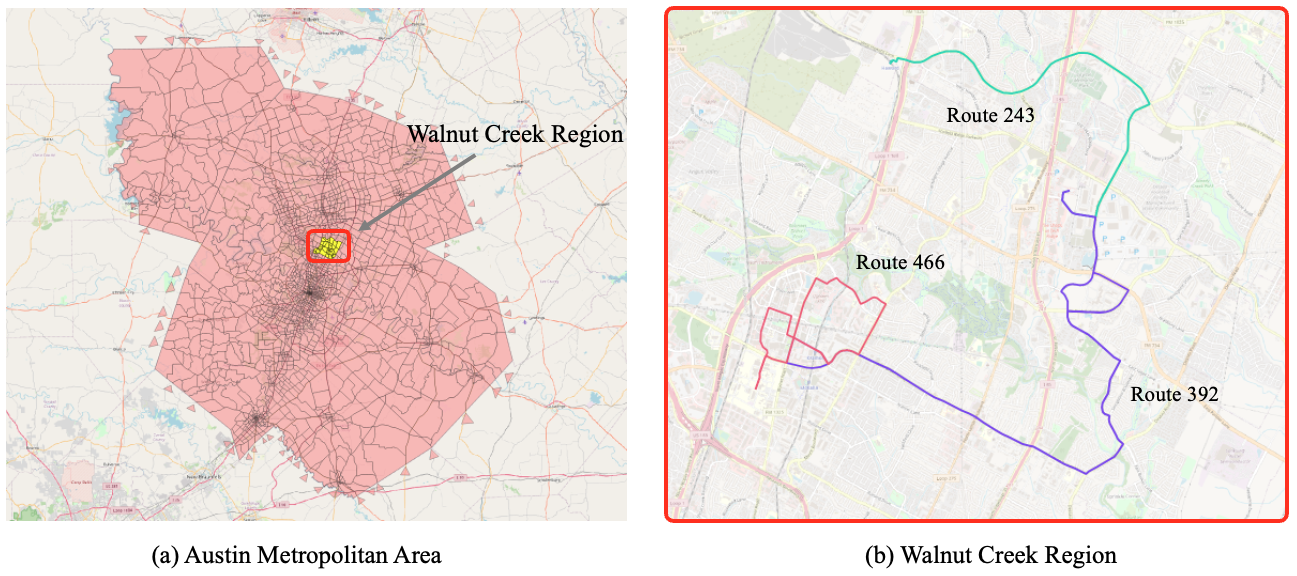}
\caption{Research area of interest}
\label{fig:exp_austin}
\end{figure}

This section assesses the efficacy of our proposed ODMTS design through a case study focused on the Walnut Creek Region located within the Austin Metropolitan Area, Texas, United States. This region is visually depicted in Figure~\ref{fig:exp_austin}. We conduct numerical experiments based on the POLARIS simulation platform, as previously introduced in Section 4, and incorporate our proposed dynamic shuttle dispatch algorithm.

To ensure the model accurately reflects inter-regional trips, we have integrated the entire Austin Metropolitan Area into the POLARIS system. A summary of our instance statistics is provided in Table~\ref{tab:instance}, including imported fixed-route transit line scheduling data from the General Transit Feed Specification (GTFS).

\begin{table}[!ht]
\caption{Instance Statistics}\label{tab:instance}
\begin{center}
\begin{tabular}{l l l l}
\hline
& Austin Metropolitan Area & Walnut Creek Region \\\hline
Transportation network & 13,177 nodes, 15,833 links & 255 nodes, 331 links \\
Traffic zone & 2,161 zones, 5,377.4 square miles & 30 zones, 13.7 square miles \\
Population & 829,962 households & 31,062 households \\
Fixed-route transit & 75 lines & 3 intra-regional lines \\\hline
\end{tabular}
\end{center}
\end{table}

Upon providing scenario inputs, POLARIS starts each simulation run with a population synthesis. This synthesized population forms the basis for the subsequent activity-based model, which generates simulated trip activities. The model subsequently utilizes trip planning and mode choice models to select suitable travel modes for fulfilling these activities. Thereafter, a traffic simulation is performed over time. Our dynamic on-demand shuttle dispatch algorithm, which is embedded within POLARIS, is triggered every 30 seconds to allocate trip requests specific to the Walnut Creek Region to the designated shuttle fleets. The CPLEX solver is employed to solve the model~\ref{formulation:dispatching_model}.

\begin{figure}[!ht]
\centering
\includegraphics[width=.8\textwidth]{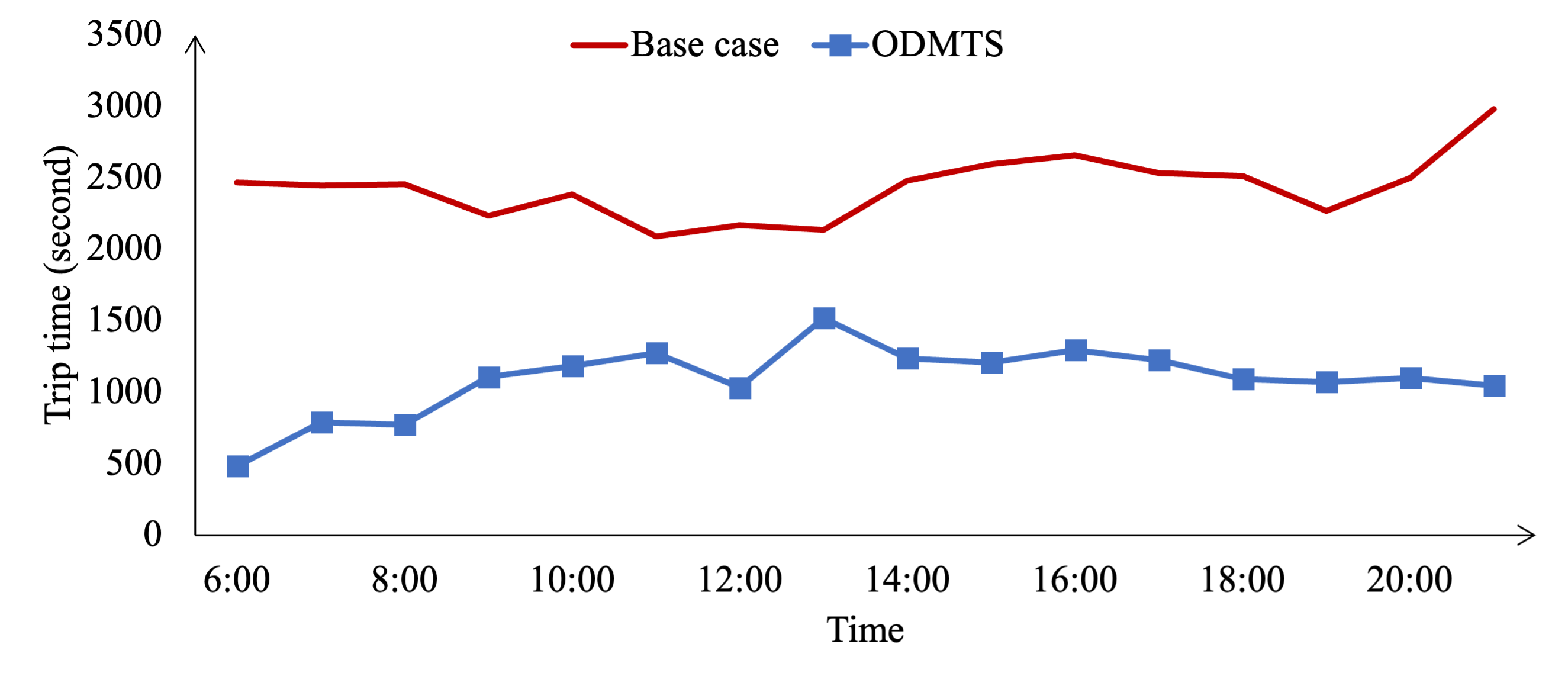}
\caption{Passenger trip travel time comparison between base case and ODMTS}
\label{fig:bus_odmts}
\end{figure}

In the context of this study, three low-ridership intra-regional bus lines (Route 243, Route 392, and Route 466) in the Walnut Creek region were removed. At the same time, we introduced a fleet of on-demand shuttles to accommodate intra-regional travel and link with the public transit systems outside of the region. Figure~\ref{fig:bus_odmts} illustrates the comparative analysis of passenger trip travel time between the base case and the ODMTS, where the ODMTS deploys five on-demand shuttles. The average passenger trip time significantly decreased from 2430.7 seconds (base case) to 998.4 seconds (ODMTS).

Under the base case scenario, transit-dependent passengers would need to walk to stations and wait for service. Conversely, in the ODMTS scenario, passengers can avail of shuttles for intra-regional trips directly to their destinations. For longer journeys, passengers can request on-demand shuttles to connect with fixed-route transit stations. It is noteworthy that the ODMTS utilizes fixed-route public transit systems to cater to concentrated corridor volumes. As a result, these systems can maintain high frequencies, thereby minimizing passenger waiting times at stations. Figure~\ref{fig:bus_odmts} reveals that passengers experience longer trip times during peak hours in the ODMTS, while they endure extended trip times during off-peak hours in the base case. The observed differences are attributed to the higher demand for on-demand services during peak hours in the ODMTS, resulting in longer waiting times. Conversely, in the base case, the reduction in travel demand during off-peak hours typically results in lower frequency fixed-route transit service to conserve operational costs, thus increasing passenger wait times at stations.

Among the removed bus lines, Route 243 and Route 392 both operated in a two-way manner with one-way journey times of 30 and 45 minutes respectively and had service frequencies of every 35 minutes. Route 466 operated in a circular fashion, with a journey time and service frequency of 35 minutes and every 30 minutes, respectively. Hence, in the previous arrangement, it required at least $\lceil 30/35 \rceil \times 2 + \lceil 45/35 \rceil \times 2 + \lceil 35/30 \rceil = 8$ buses to serve the three lines. However, the on-demand ride system necessitates only five shuttles. If we consider utilizing standard-sized buses as shuttles in ODMTS, operational costs would see a reduction of $(8-5)/8 \times 100=37.5$ percent. It's worth noting that, in a real-world context, fixed-route buses are typically larger than what would be required for shuttles, suggesting that the estimated 37.5 percent is a conservative approximation of potential operational cost savings. This finding underscores the substantial advantages of implementing an on-demand service under the ODMTS, particularly in regions with low ridership.

\begin{figure}[!ht]
\centering
\includegraphics[width=.8\textwidth]{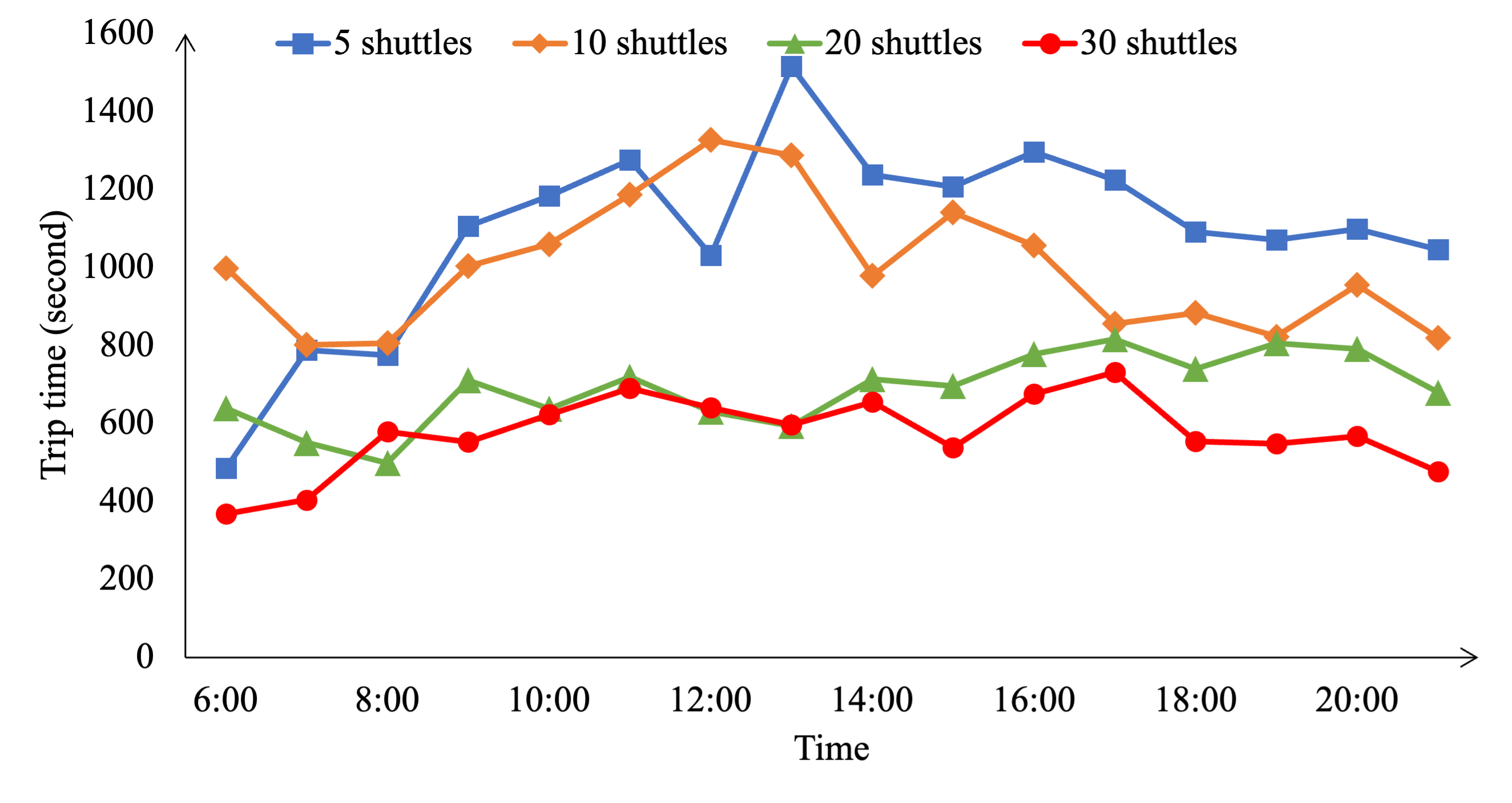}
\caption{Passenger trip travel time with different number of shuttles}
\label{fig:num_of_shuttles}
\end{figure}

Figure~\ref{fig:num_of_shuttles} elucidates how passenger trip travel times alter with varying numbers of shuttles. As the quantity of shuttles increases, passenger trip times decrease. During peak hours, scenarios with 5 and 10 shuttles witness longer passenger trip times due to a shortage of shuttles. Conversely, scenarios with 20 and 30 shuttles display more stable trip times. These findings provide valuable insights for transportation planners aiming to strike a balance between operational costs and passenger trip experience.

%% file: sections/6_conclusions.tex
\section{Conclusions}
\label{sect:conclusions}

This study focuses on the design and operation of public transit systems within low-ridership regions, particularly exploring the promise of On-Demand Multimodal Transit Systems (ODMTS) in enhancing accessibility and convenience for public transit users. The ODMTS framework proposed advocates for a strategic transformation: substitute local fixed-route transit lines within low-ridership areas with a dynamic fleet of on-demand shuttles. The integrated system operates utilizing fixed-route transit services to cater to concentrated corridor volumes, while deploying on-demand ride services to bridge the first-mile and last-mile gaps in travel. This enables each mode of transit to play to its strengths. The study presents an in-depth discussion on the overall system design and the dynamic methodologies employed in shuttle dispatching. To evaluate the efficacy of the proposed model, we conducted numerical experiments within the Austin Metropolitan Area, Texas, United States, utilizing the POLARIS simulation platform with the proposed shuttle dispatching algorithm embedded. The results showed that, compared to the current system, the proposed ODMTS can significantly enhance the public transit travel experience, while reducing operating costs. The insights derived from this study support designing and operating more flexible, accessible, and convenient public transit systems in the upcoming era of Mobility-as-a-Service.

%% file: sections/acknowledgements.tex
\section{Acknowledgments}
This article and the work described were sponsored by the U.S. Department of Energy (DOE) Vehicle Technologies Office (VTO) under the Systems and Modeling for Accelerated Research in Transportation (SMART) Mobility Laboratory Consortium, an initiative of the Energy Efficient Mobility Systems (EEMS) Program. The U.S. Government retains for itself, and others acting on its behalf, a paid-up nonexclusive, irrevocable worldwide license in said article to reproduce, prepare derivative works, distribute copies to the public, and perform publicly and display publicly, by or on behalf of the Government.